\documentclass[12pt, letterpaper]{article}

\usepackage{wasysym}
\usepackage{booktabs}
\usepackage[T1]{fontenc}
\usepackage{graphicx}
\usepackage[latin,english]{babel}
\usepackage{ledmac,ledpar}
\lineation{section}
\usepackage{graphicx,vmargin}
\setpapersize{USletter}

\newcommand{\exces}{\rotatebox{90}{\dag}}
\setmarginsrb{0.8in}{0.8in}{0.8in}{0.9in}{0pt}{0pt}{0pt}{24pt}

\begin{document}

\begin{center}
\textbf{\begin{Large}\textit{Solutio Theorematum} \\by Adam Adamandy Kocha\'nski -- Latin text with annotated English translation
\end{Large}}
\end{center}
\begin{center}
translated by Henryk Fuk\'s\\
Department of Mathematics,
Brock University, St. Catharines, ON, Canada\\
email \texttt{hfuks@brocku.ca}\end{center}

\noindent \textit{Translator's note:} The Latin text of \textit{Solutio theorematum} presented here closely  follows the original text published in 
\textit{Acta Eruditorum} 
\cite{KochanskiSolutio}. Punctuation, capitalization, and mathematical notation have been preserved. Several misprints which appeared in the original are also reproduced unchanged, but with a footnote indicating correction. Every effort has been made to preserve the layout of original tables.
The translation is as faithful as possible, often literal, and it is mainly intended to be of help to those who wish to study the
original  Latin text. In the appendix, all propositions from Euclid's \emph{Elements} mentioned in the text are listed in both Latin
and English version.

\setlength{\Lcolwidth}{0.465\textwidth}
\setlength{\Rcolwidth}{0.495\textwidth}

\begin{pairs}
\begin{Leftside}
\beginnumbering
\selectlanguage{latin}
\noindent\setline{1}
\pstart
\begin{center}
 \emph{\Large SOLUTIO THEOREMATUM Ab illustri Viro in Actis hujus Anni Mense Januario,
pag. 28. propositorum, data ab Adamo Adamando Kochanski S.J. quondam Pragensi
Mathematico.}
\end{center}

DUPLICATIONEM Trigoni Isogoni, citra Proportiones demonstrandam P.
Sigismundus Hartman\footnote{Sigismundus Ferdinandus Hartmann SJ (1632--1681)
-- Bohemian Jesuit and mathematician, professsor of the University of Prague.} e Soc. Jesu, publico Programmate proposuerat, istudque
Schediasma mihi pro veteri necessitudine transmiserat e Bohemia in Poloniam.
Reposui humanissimo Authori solutiones fere vicenas, ab eo sic probatas, ut una
cum aliis, aliunde ad se missis, in lucem daturus fuisset, si Parcae Viro tanto
pepercissent.

Cum vero his primum diebus in manus meas venerint Acta Eruditorum Lipsiensia, \&
in his Solutio problematis istius Hartmanniani, a quodam illustri viro data, \&
ad P. Coppilium\footnote{Matthaeus Coppilius SJ (1642--1682) -- Bohemian Jesuit and mathematician,
author of books on mechanics.}, defuncti in Mathematici munere successorem, quodammodo
directa, quasi is editioni posthumae operum Hartmanni, maximeque \emph{Protei
Geometrici}, ab eo nuper promissi incumberet ; cum tamen ab obitu Authoris
elaboratum nihil, sed prima solum Operis lineamenta reperta fuisse mihi ab
Amicis nunciatum fuerit ; Eam ob rem non aegre laturum spero P. Coppilium, si
dum illum alio in opere fructuose versari intelligo, hic ejus partes occupare
ausus fuero: non enim it temere, aut praesidenter egisse videbor, sed veluti
jure quodam antiquitatis ; quod videlicet ante illum in eodem Matheseos Pragensis
pulvere quendam Professor fuerim versatus; adeoque prior tempore, licet non
eruditione.

Dilatis autem in aliud tempus meis illis Duplicationum Particularium, \&
Universalium Solutionibus, una cum Pythagoricae nova, ac multiplici
Demonstratione, ceterisque meis considerationibus Geometricis, usque dum
ultroneam Typographi, vel Bibliopolae cujuspiam humanitatem invenerint;
suffecerit hoc loco strictim ea persequi, quae pertinent ad geminum illud
Theorema Geometricum, quod Anonymus ille Problematis Hartmanniani 
$\delta\epsilon \acute{\iota} \chi\tau\eta\varsigma$
 loco
sup. memorato proposuit.
\pend
\pstart

\begin{center}\textbf{\Large ARTICULUS I.}\end{center}

Circa primum illorum Theorematum consideranda veniunt sequentia. Inprimis
dissentire me in eo ab \textit{Illustri Viro}, quod is existimet Theorema Pythagoreum
continuari in Sphaera, Pyramidem Rectangulam circumscribente; cum nec Pythagoras
suum illud Orthogonium, tanquam circulo inscriptum consideraverit, nec, si
universaliter agamus de potentiis cujusvis Trianguli, haec consideratio proprie
ad Circulum pertinere videatur, sed potius ad Parallelogrammum universim ;
quando videlicet istud Diametro sua sectum concipitur in duo Triangula aequalia,
eaque Orthogonia, Amblygonia, vel Oxygonia, pro diversitate parallelogrammi;
Harum enim Diametrorum Potentiae cum suis lateribus comparantur 47. \textit{Primi}, nec
non 13, \& 13. \textit{Secundi Elem.} Universalius autem hoc ipsum consideratur 31 \textit{Sexti}
saltem quoad orthogonium Triangulum ; nam quoad reliqua, videri possunt ea, quae
demonstrat Clavius e Pappo, \textit{in Scholio ad 47. Primi, ad finem.}

Quamobrem non ineleganter fluente Analogia, nobis dicendum videtur, Pyramides
Triangulares Orthogonias, Amblygonias, \& Oxygonias, cum suis laterum, ac
diametrorum Potentiis, immediate quidem reduci oportere ad Prismata sua, bases
triangulares habentis, quorum Pyramides illae sunt partes tertiae \textit{per Prop. 7.
Duodecimi Elem.} Haec ipsa vero Prismata, tanquam partes revocantur ad totum
Parallelepipedum, cujus sunt medietates. Istas porro relationes partium ad sua
Tota intelligi volumus de ordine \textit{Doctrinae} potius, quam \textit{Naturae},
constat enim, Triangulum esse prius, ac simplicius Parallelogrammo, non minus
ac Pyramidem Tetrahedram Prismate, vel Parallelepipedo: Hinc 
$ \delta o \gamma \mu \alpha \tau \iota \chi \tilde{\omega} \varsigma$
 sperando,
\& ipsae Linearum Potentiae, non Triangulis aequilateris, sed Quadratis, omnium
consensu taxantur, licet illa sint istis priora, magisque simplicia.

Quamvis autem Pyramidum Polyhedrarum aliae inscribi possint Sphaeris , aliae
vero Sphaeroidibus, ac tum Potentiae laterum conferri cum Diametro corporis
circumscribentis: eadem tamen Pyramides adhuc secari poterunt in Tetrahedras,
atque ita Prismatibus suis, ac tum Parallelepipedis veluti postliminio quodam
restitui. Si quis nihilominus omnia Trilatera a Quadrilateris, omnia Tetrahedra
a Pentahedris \& Hexahedris emancipare  contenderit cum eo nequaquam cruento
Marte dimicabimus.
\pend
\pstart

\begin{center}\emph{\Large PROPOSITIO I. Theorema.}\end{center}
{\large In omni Pyramide rectangula, tria Quadrata laterum, Angulum rectum in vertice
comprehendentium, aequalia sunt Quadrato Diametri totius Parallelepipedi aeque
alti, Pyramidem illam complectentis.}

Sit Pyramis ABC rectangula ad verticem D. sive jam sit \textit{Aequilatera}, prout in
Cubo \textit{Figura} I. sive \textit{Isosceles}, ut est in Fig. II Parallelepipedi, supra basin
quadratam ADCF assurgentis ; sive demum \textit{Scalena}, qualem exhibet in \textit{Fig.} III.
solidum rectangulum, supra basin ADCF altera parte longiorem, erectum.

Dico, tria Quadrata DA. DB, DC aequari Quadrato Diametri AE, per oppositos
solidi angulos incedentis. Nam primam in Triangulo ADB angulus D rectus est, ex
hypoth. Igitur Qu. lateris AB aequatur (\textit{47. I. Elem.}) Quadratis AD. DB duorum
laterum datae Pyramidis. Deinde Triangulum pariter ABE rectangulum est ad B. (id
ostendi potest \textit{per 4. Undecimi}) Quocirca Quadratum AE aequabitur Quadrato AB,
hoc est duobus DA. DB, \& insuper Quadrato BE, hoc est ipsi aequali DC, quod est
tertium latus datae Pyramidis ABCD.

Tria igitur omnis Pyramidis rectangulae latera, potentia aequantur Diametro
Parallelepipedi Pyramidem continentis, q. e. d.
\pend
\pstart
\begin{center}\emph{\Large Corollarium I.}\end{center}
 Colligitur hinc, in omni Prismate rectangulo ACDBGA, cui basis
est Orthogonium Triangulum ADC cujus Parallelogrammi rectanguli AGEC, quod
angulos D. \& B. rectos subtendit, Diametrum AE. Potentia aequari iisdem tribus
lateribus Pyramidis rectangulae ABCD: eo quod ipsa DE sit eadem omnino cum
Diametro totius Solidi GC.
\pend
\pstart
\begin{center}\emph{\Large Corollarium II.}\end{center}
Hinc quoque manifestum est, Diametrum Sphaerae, quae Pyramidi rectanguli ABCD
circumscripta est, aequari potentia tribus lateribus ejusdem Pyramidis, rectos
angulos constituentibus in vertice D. Nam cum omni solido rectangulo Sphaera
circumscribi possit, non minus ac Circulus ejus basi rectangulae, erit Diameter
solidi eadem omnino quae circumscriptae Sphaerae, per ea, quae Pappus demonstrat
\textit{Lemmate 4. apud Clavium ad calcem Lib. 16. Elem.}

Si cui placuerit in simili materia ingenium exercere circa Pyramides
Triangulares, tam Acutangulas, quam Obtusangulas, itemque mixtis in vertice
angulis contentas  in illis Potentias laterum cum Diametro totius Solidi
obliquanguli comparando; id vero non difficulter poterit expedire ope duarum
Propositionum, videlicet penultimae, \& antepenultimae \textit{Lib. Secundi Elem.}

\pend
\pstart
\begin{center}\textbf{\Large ARTICULUS II.}\end{center}
Circa alterum Theorema percontatur \textit{Illustris Vir} AN sicut in Circulo unica Media
Proportionalis, ita etiam insistendo Analogiae, in Sphaera duae Mediae inveniri
possint ? Ad hanc quaestionem Respondeo I. Nec Antecedens Analogiae hujus tam
ratum esse, ut absolute loquendo, Circulus duas Medias excludere pronunciari
possit, aut debeat: Fieri namque potest, ut in Semicirculo ACD (\textit{inspice Fig. 4.})
continuae sint AB. BC. BD. DA. cujus Problematis Geometricam constructionem
Mathematum peritis propono, interim vero Arithmeticum sequentibus numeris
expono. Ponatur enim
\pend
\end{Leftside}

\begin{Rightside}
\beginnumbering
\selectlanguage{english}
\noindent\setline{1}
\pstart
\begin{center}
 \emph{\Large SOLUTION OF THE THEOREM proposed by an illustrous man in this year's January issue  of Acta
on p. 28, given by Adam Adamandy Kochanski SJ, once mathematician of Prague.}
\end{center}

Fr. Sigismundus Hartman from Soc. of Jesus proposed in a public program
the problem of DUPLICATION of an equilateral triangle without the use of proportions,
and conveyed this question to me from Bohemia to Poland, on the account of old friendship.
I returned to the kindest author perhaps twenty solutions, by him thus examined,
so that one of these, together with others, send to him from elsewhere, was to be
put to light, if only Fate had spared the man so great.

When, however, those days  Acta Eruditorum of Leipzig came into my hands, 
and in them solution of this problem of Hartman, given by a certain illustrous man,  
directed to Fr. Copillus, successor of the  deceased in the mathematical office,
and arranged in a certain way, as if he was leaning towards posthumous edition
of Hartman's works, especially \textit{Protei Geometrici}, recently promised by him;
when still from the death of the Author nothing has been worked out,
but friends announced to me that only the preliminary outline of the Works has been
produced; therefore I hope that Fr. Copillus is not offended if, while
he, as I understand, is busy with other fruitful works, I dared
to assume his role: this will not, indeed, seem to be acting rashly, or daringly,
but  by a certain old right; because evidently I used to dwell in the same
dust of Prague as  Professor of Mathematics\footnote{Kocha\'nski
stayed in Prague from 1670 to 1672, lecturing in mathematics 
and (probably) moral philosophy.}; indeed, preceding [Fr.
Hartman] in time, but not in erudition. 

Leaving for another time my particular and general solutions of the Duplication,
one with a novel Pythagorean proof, others based on  my Geometrical considerations,
until they find a willing printer or a kind publisher; it will suffice in this 
place to pursue ony what pertains to these two Geometrical Theorems, which
the Anonymus exhibitor put forward in the aforementioned place. 
\pend
\pstart

\begin{center}\textbf{\Large ARTICLE I.}\end{center}

The following considerations revolve around the first of these theorems.
First of all I disagree with the \emph{Illustrious Man} in his opinion that the Pythagorean 
Theorem extends to a sphere circumscribed around a right-angled pyramid; when neither
Pythagoras thought of his right triangle as if it was inscribed in a circle, nor,
if we talk in general about properties of an arbitrary triangle, does
such considerations seem to pertain specifically to a circle, but rather
to parallelogram in general; when clearly this [parallelogram] is understood
to be cut by its diameter [i.e., diagonal] into two equal triangles, either right-angled or
obtuse-angled or acute-angled, according to diversity of parallelograms;
Powers [squares] of these diameters [diagonals] are compared with powers of their
[parallelograms] sides in prop. 47 of the \emph{first book of Elements} as well as
in  and prop. 13 and prop. 13\footnote{Obviously, this should be ``prop. 12. and prop. 13''.} 
of the \textit{second book}.
This is considered more generally in prop. 31 of the \textit{sixth book}, as long as the right
triangle is considered. On the other hand, cases which remains 
can be seen as those demonstrated by Clavius following Pappus \textit{at the
end of commentary to prop. 47}.\footnote{Pappus builds parallelograms
on sides of an arbitrary triangle, cf. p. 366 of \cite{Euclid}.}

It seems we need to state an elegantly flowing analogy, that it is right
to put triangular prisms, right-angled as well as obtuse-angled and acute-angled,
with powers [squares] of their sides and diameters, into their respective prisms,
having triangular bases, of which these prisms constitute third parts [by volume],
by \textit{prop. 7 of the twelfth book}. Such prism are in truth recalled just as if they were
parts of a complete parallelepiped, of which they are halves. Hereafter we want these
relations of parts to their wholes  to be understood from
\textit{Principles} rather  than from \textit{Nature}, as it is evident that the triangle is prior to and
simpler than a parallelogram, not less as the tetrahedral pyramid is prior to and simpler
than a prism, or parallelepiped: from there observing the thing dogmatically, 
powers of own lines taken together are valued not with equilateral triangles, 
but squares, although  those [triangles] are prior to these [squares], and much simpler.

Although some polyhedral pyramids could be inscribed in spheres, and
  others in spheroids [i.e., ellipsoids of revolution], and then powers of their sides
could be matched against diameters of circumscribed bodies: besides, the same pyramids 
could be divided into tetrahedrons, and therefore also brought into their
own prisms, and then parallelepipeds, as if by the right to return home.
None the less, if somebody wanted to alienate all trilaterals from quadrilaterals,
all tetrahedra from pentahedra and hexahedra , we will by no means 
spill blood in a fight with him.
\pend
\pstart

\begin{center}\emph{\Large PROPOSITION I. Theorem.}\end{center}
{\large
In every right-angled pyramid, sum of squares of three sides coming from
the vertex embracing the right angle is equal to the square of the diameter 
of the complete parallelepiped of equal height, encompassing this pyramid.} 

Let the pyramid ABC be right-angled at the vertex D, whether \textit{equilateral},
as in the cube of Fig. I, or   \textit{isosceles}, as in the parallelepiped of
Fig. II, rising over the square basis ADCF, or finally \textit{sceles}\footnote{Having
three equal sides.}, as in the rectangular solid erected over the basin
ADCF elongated in the other direction, as displayed in Fig. III.

I say that three squares of DA, DB, and DC are equal in sum to the square
of the diameter AE, stretched between opposite angles of the solid.
For, first of all, in the triangle ADB, the angle D is right, by hypothesis.
Therefore the square of AB is equal (by \textit{prop. 47 of Elem. I}) to the sum of squares
of two sides AD and DB of the given pyramid. Next also triangle ABE is
right-angled at B (this can be shown by \textit{prop. 4 of the eleventh book}), on account of
which square of AE is equal to the square of AB, that is sum of squares
of DA and DB, plus square of BE, which itself is equal to DC, the third
side of the given pyramid ABCD.\footnote{$AB^2=AD^2+DB^2$, $AE^2=
AB^2+BE^2=AD^2+DB^2+BE^2=AD^2+DB^2+DC^2$.}

Therefore three sides of any rectangular pyramid are equal
in power\footnote{i.e., sum of their square is equal to the square of...}
to the diameter of the parallelepiped enclosing the pyramid, Q.E.D.
\pend
\pstart
\begin{center}\emph{\Large Corollary I.}\end{center}
One obtains from there that in all rectangular prisms ACDBGA, whose
base is a right triangle ADC, square of the diameter AE of the right-angled
parallelogram AGEC, which extends below right angles D and B, is equal to
sum of squares of three sides of the right-angled pyramid ABCD:
consequently DE\footnote{Clearly a misprint. Should be $AE$.} itself would be entirely the same as the diameter
GC of the whole solid. 
\pend
\pstart
\begin{center}\emph{\Large Corollary II.}\end{center}
From this it is evident that the square of the diameter of the sphere circumscribed
around a right-angled pyramid ABCD, is equal to the sum of squares
of three sides of this pyramid forming the right angle at the vertex D.
For while every rectangular solid can be circumscribed by a sphere,
as much as its rectangular base [can be circumscribed] by a circle,
the diameter of this solid will be entirely the same as the diameter
of the circumscribing sphere, as  Pappus demonstrated in
\textit{Lemma 4 in Clavius' commentary at the end of book 16\footnote{Book XVI was a medieval addenum to \emph{Elements}.}
 of Elements.}

If somebody would please to exercise his talent in similar matters
regarding triangular pyramids, either acute-angled or obtuse-angled,
and likewise stretching over mixed angles, comparing sums of squares
of the sides with the diameter of the whole solid; this will
certainly not be difficult to obtain by the power of two propositions,
namely the second last and the third from the end proposition of
the \textit{second book of Elements.}

\pend
\pstart
\begin{center}\textbf{\Large ARTICLE II.}\end{center}
The Anonymous \textit{Illustrious Man} inquires about another theorem,
as if there exist one geometric mean in a circle, could two
means me found in a sphere? I make first response to this question,
that the first part of this analogy is not so strongly established that,
absolutely speaking,  it would exclude a possibility that
a circle  having two means.
For in fact, it can happen that in a semicircle ACD (\textit{see Fig. 4})
there are successive lines AB, BC, BD, DA. I leave Geometric Construction
of this problem to mathematical experts, and in the meanwhile
I explain  Arithmetic one with the following numbers. 
It is namely assumed that diameter will be
\pend
\end{Rightside}

\Columns

\end{pairs}

\begin{center}
\begin{tabular}{clllrlr}
Diameter & AD. & - & - & 2 00000 & 00000 &  \\
erit & AB. & - & - & 63534 & 43923. & \exces\\
 & BC. & - & - & 93114 & 24637. & \exces\\
 & BD. & - & - & 1 36465 & 56077. & \exces
\end{tabular}
\end{center}

\begin{pairs}
\begin{Leftside}
\selectlanguage{latin}
\pstart
Unde \textit{per 19. Septimi Elem.} erunt aequalia Rectangula
\pend
\end{Leftside}

\begin{Rightside}
\selectlanguage{english}
\pstart
From there by \textit{Prop. 19 of the seventh book of Elements} [the following] rectangles will be equal
\pend
\end{Rightside}

\Columns
\end{pairs}

\begin{center}
\begin{tabular}{clllrlr}
DAB. &- &- &1 17068 & 87846 &00000 &00000.\\
CBD. &- &- &1 17068 & 87846 &55798 &69049.
\end{tabular}
\end{center}
\begin{center}Gg.\end{center}	

\begin{pairs}
\begin{Leftside}
\selectlanguage{latin}
\pstart
Et \textit{per 20. ejusdem}, Quadrata mediarum aequabuntur Rectanguli sub earundem
extremis.
\pend
\end{Leftside}

\begin{Rightside}
\selectlanguage{english}
\pstart
And \textit{by Prop. 20}\footnote{Prop. 20 is often omitted in modern editions of \emph{Elements}, as it is considered 
a later addition, and a direct consequence of proposition 19.}
 \textit{of the same}, squares of means will be equal to [areas of] rectangles
under their outer segments.\footnote{``Outer'' means neighbours in the sequence. Kochanski considers sequence of linear segments
AB, BC, BD, DA, where both BC and BD are geometric means of their nearest neighbours in the sequence.}
\pend
\end{Rightside}

\Columns
\end{pairs}

\begin{center}
\begin{tabular}{llllrlll}
 $\Square$BC &-&-&-& 86702 &62877 &05305 &81769.\\
        ABD &-&-&-& 86702 &62877 &72943 &70071.\\
  $\Square$BD &-&-&-& 1 86288& 49276& 27056& 29929.\\
       ADBC &-&-&-& 1 86228& 49274& 00000& 00000.
\end{tabular}
\end{center}

\begin{pairs}
\begin{Leftside}
\selectlanguage{latin}
\pstart
Respondeo II.
 De Analoga illa nihil certi statui posse videtur: Unius enim
Mediae inventio, quae Circulo tribuitur, etiam Sphaerae congruit ; \& inventio
duarum Mediarum, hic a nobis demonstranda, aeque ad circulum, sicut \& Sphaeram
aptari poterit, quemadmodum ex dicendis constabit.
\pend

\pstart
\begin{center}\emph{\Large PROPOSITIO II. Theorema.} \end{center}
In Circulo ADC a Diametro AC descripto ducantur utcunque ad peripheriam duae
rectae AD DC: tum ex D cadat ad AC perpendicularis D; similiterque ex E sit ad
AD normalis EF.

Dico in Circulo ADC, haberi Quatuor continue Proportionales, AF. AE. AD. AC.
Describatur enim Diametro AD Circulus AGDE, quem EF producta secet in G, \&
connectantur G A. GD.

\emph{Demonstr.} Recta AD subtendit Angulum rectum DEA. Igitur Circulus AGDE Diametro
AD descripsit transit per verticem Anguli recti DEA. \textit{juxta  Schol. Clavii ad 31.
3 Elem.} Est autem recta EFG perpendicularis ad Diametrum AD. Ergo \textit{per 3. 3.
Elem.} tota EG bifariam secatur in F. Triangula igitur AFG, AFE orthogonia in F,
sunt \textit{per 4.1. Elem}. invicem aequalia : Eademque de causa aequantur Triangula
DFG. DFE, ac proinde \& totum DGA toti DEA aequale. Jam sic. In Orthogonio AED
(par ratio de aequali AGD) ab angulo recto E cadit perpendicularis EF in basin
AD: Ergo \textit{per Coroll. 8.6. Elem.} Proportionales sunt tres AF.AE.AD. Sed eadem de
causa in Orthogonio ADC duabus postremis e praecedenti serie, videlicet AE. AD.
proportionalis est tertia AC. Igitur omnes quatuor AF.AE.AD.AC sunt in continua
Proportione intra Circulum ADC, q.e.d.
\pend

\pstart
\begin{center}\emph{\Large Corollarium.} \end{center}

Non difficulter hinc elicitur, easdem quatuor continuas in Spaerae  quoque
concipi posse, non modo praedicta, sed \& alia ratione: Fingamus enim Sphaera ADC
auferri Segmentum AHDA, cujus basis erit
Circulus a Diametro AD, cujus meditas esto AGDA. Ductis autem
Orthogonalibus DE. EF. in plano Circuli Spherae maximi ADC, nec
non Orthogonali FG, in altero plano Semicirculi AGD, hoc est basi 
Segmati AHDA, jungatur AG: erunt enim ut antea, quatuor AF.
AF, AD.AC. continue proportionales, id patet e praecedenti discursu,
qui non difficulter huc applicari poterit, licet plana Circulorum ADG.
AGD. sint diversa, \& ad rectos invicem collocate.

Notandum vero est, Theorema praecedens, loquendo pressius, non tam Semicirculo,
quam Orthogonio cuilibet in similia subdiviso convenire: quia tamen ejus
Demonstratio sequentibus inserviet Problematibus, visum est illud hoc loco
tantisper indulgere Circulo, sine quo illa absolvi non possunt.
\pend

\pstart
\begin{center}\emph{\Large PROPOSITIO III. Problema.} \end{center}

{\Large Inter duas datas, duas medias in continua ratione, duobus tantum digitis
reperire.}

Celeberrimum illud Problema Deliacum quot \& quanta totius Orbis eruditi
exercuerit ingenia, Geometris est notissimum, ut \& variae  illius absolvendi
Praxes Organicae, a compluribus excogitatae, quarum aliae aliis sunt
operosiores: Nostra haec videri poterit nonnihil Paradoxa, quod duobus tantum
digitis unius manus, absolvatur, cum nonnullae requirant, \& occupent utramque.

Datae sint, \textit{in Figura VI.} duae AC. AB. quas inter duae mediae quaeruntur. In
communi utriusque termino A figatur Regula AZ, instructa Cursore FY, qui semper
insistat ad rectos ipsi regulae AZ, idque firmiter, ubicunque collocetur. In
illa fumatur AF, qualis datae AB, minori altera AC. Descripto autem super tota
AC Semicirculo ADC in plano quopiam verticali, hoc est ad Horizontem recto, cui
aequidistet Diameter AC: applicetur ad Peripheriam ADC Stylus quidam gracilis
DS, e quo deorsum propendeat filum subtile cum appenso Pondere X, vel certe
hujus loco regula quaedam sub gravis, accurate tamen aequilibrata: Nam si Stylus
DS duobus digitis apprehensus pedetentim promoveatur per Circumferentiam AD,
usque dum Perpendiculum DE cum Cursore FE sese mutuo intersecent alicubi in
recta AC, velut in Puncto E: istud probe notatum offeret quatuor Proportionales,
quarum duae AE. AD. inter datas extremas AF hoc est AB nec non AC interponentur.

Demonstratio Problematis hujus, quoad rem, eadem est; cum adducto praecedenti
Theoremate.
\pend

\pstart
\begin{center}\emph{\Large PROPOSITIO IV. Problema.} \end{center}
\pend

\pstart
\begin{center}{\Large Id ipsum aliter, una Circini apertura}.\end{center}

Quoniam praecedens praxis ob situm plani verticalem, \& usum Perpendiculi,
nonnihil impedita cuipiam videri possit, quin \& a Geometriae moribus aliena,
dabimus alteram, praedictis incommodis haud obnoxiam.

Positis iisdem, in locum perpendiculi DEX Figurae praecedentis, subrogetur \textit{in
hac Figura VII.} Parallelogrammum materiale KLMN, cujus unum Latus KL. fixum sit
in plano, alterum vero MN mobile, semper tamen ad rectos ipsi AC. Nam si tota AC
divisa bifariam in O Circini pes unus figatur in O, alter autem intervallo OC
diductus, tamdiu in Arcu CD provehatur, impellatque binas Regulas AZ. \& MN in
puncto intersectionis D, usque dum regula MN Cursorem FY secet in puncto E,
posito in recta  AC, obtinebuntur eaedem mediae AE. AD. longe commodiori
ratione,  quam fuerit praecedens;  quae tamen ipsa, si Figura magnae molis
fuerit, in vasto quopiam pariete  usui esse poterit Architectis.

Non est cur hoc loco moneam de circino, ejusque Cruribus in regula quapiam
mobilibus, nec de acie pedis alterius, quae regulas in puncto D subtiliter
impellere debet; nec denique de nisu quodam regularum AZ, MN contra Circinum ;
hunc enim vel ipsarum regularum  pondere, vel Elatare quopiam, aut unius digiti
impulsu consequetur ingeniosus quivis, ac istis etiam nostris longe meliora
excogitabit.
\pend
\end{Leftside}

\begin{Rightside}
\selectlanguage{english}
\pstart
Secondly, I respond than nothing certain seems to possible to state about this analogy.
Invention of one mean, assigned to a circle, suits to sphere too.
Invention of two means, demonstrated here by us, can be adapted as well to
a circle as to a sphere, in what way it will agree with what has been said.    
\pend

\pstart
\begin{center}\emph{\Large PROPOSITION II. Theorem.} \end{center}
In a cirlce ADC traced out from a diameter AC two straight lines are drawn
as far as to the perimeter: then from D a line D\footnote{A clear misprint: this should be DE.} perpendicular to AC line is led,
and similarly from E a line EF normal to AD.

I say that in the circle ACD four consecutive proportionals exist, AF, AE, AD, AC. Indeed, let
a circle AGDE with diameter AD be drawn, which intersects with extension of the line EF
at G, and let G connecting lines GA and GD be made.

\emph{Proof.} Line AD extends beneath the right angle DEA. Therefore the circle AGDE determined
by the diameter AD passes through the vertex of the right angle DEA, \textit{according to
commentary of Clavius to prop. 31 of book 3 of Elements}. The straight line EFG
is in fact perpendicular to the diameter AD. Therefore by \textit{prop. 3, book 3 of 
Elements}, the entire EG is divided into two equal parts at F. The triangles AFG, AFE
having straight angle at F, are by \textit{Prop. 4.1 Elem.} equal to each other. For the same
reason triangles DFG, DFE are equal, and hence the whole [triangle] DGA is equal
to the whole [triangle] DEA. Now in the right-angled triangle AED (similarly reasoning
applies to AGD) from the right angle E a perpendicular line EF falls onto the 
base AD: therefore, by \textit{Coroll. 8. 6. Elem.}, there are three proportionals AF, AE, AD.
But for the same reason, in the right-angled triangle ADC, AC is proportional 
to he two last lines from the aforementioned sequence [of proportionals], namely AE and AD. 
Therefore all four AF, AE, AD, AC are proportional in succession within the circle
ADC, Q.E.D.
\pend

\pstart 
\begin{center}\emph{\Large  Corollary.} \end{center}
Form there it is not difficult to elicit that the same four successive proportionals can be
devised in a sphere, not by the preceding method, by by a different reasoning: let
us namely imagine that a segment\footnote{spherical cup} AHDA is taken from  a sphere ADC, whose basis is
a circle with diameter AD, and whose one half is AGDA. Drawing perpendicular lines
DE, EF in the plane of the great circle ADC, and perpendicular line GF, in another 
planar semicircle AGD, that is, in the base of the segment\footnotemark[\value{footnote}]
AHDA, let they be joined by AG: there will be, as before, four successive proportionals
AF, AG, AD, AC. It stands clear from the previous discourse, which could be applied here
without difficulty, that it is permitted that planes of circles ADG and AGD are different,
and  placed at right angles to each other.

One must observe, however, that the previous theorem, speaking more precisely,
is not as tied to a semicircle as to an arbitrary triangle subdivided
into similar ones: yet because its proof lends itself to the following problems,
certain leniency toward the circle  appeared in this place, without which
they could to be brought out. 
\pend

\pstart
\begin{center}\emph{\Large Proposition III. Problem.} \end{center}

{\Large Between two given [quantities], find two means in successive proportions, with only two fingers.}

It is very well known to Geometers how many and how great learned [men] 
exercised [their] talents on this most famous Delian problem,
and how  various
practical methods of its solution utilizing mechanical instruments, devised by many, 
have advantage one over another. With our method one could  these as  paradoxes,
because it utilizes only two fingers of one hand, while some other [methods]
require and occupy both [hands].

Given are, as in \textit{Figure VI}, two [quantities] AC and AB, between which two means
are sought.  At their common end A a ruler AZ is fixed,equipped with cursor FY,
which always firmly stands at the right angle to the ruler AZ, wherever placed.	
With the ruler the distance AF is taken, equal to the given AB, smaller that AC.
Somewhere in a vertical plane, that is, perpendicular to the horizontal line, 
let a semicircle ADC be drawn over the entire AC, whose diameter is equal to AC. 
Let at the circumference of ADC a thin stylus DS be placed, from which a fine string is
hanging down with a weight X attached, and at which finally the ruler is placed
in equilibrium under gravity. For if the stylus DS is carefully held with two fingers 
and moved through the circumference AD, all the way until the perpendicular DE
intersects with he cursor FE somewhere on the straight line AC, for instance at point E,
this, as it rightly to be noted, produces four proportionals, of which two, AE and AD,
are interposed between two given outward  [lines] AF (that is AB) and AC.

Demonstration of [correctness of solution of] this problem follows from the preceding proved
theorem.
\pend

\pstart
\begin{center}\emph{\Large PROPOSITION VI. Problem.} \end{center}
\pend

\pstart
\begin{center}{\Large The same thing differently, with one aperture of the compass.}\end{center}

Because the previous method seems to be hindered by the positioning of the vertical plane
and by the use of the plumbline, which is alien to customs of geometers, we will
give another [method], not burdened with the aforementioned inconveniences. 

Setting up things as before, in place of the perpendicular DEX of the preceding figure,
let in \textit{Figure VII} a material parallelogram KLMN be substituted, whose one side KL is fixed 
in the plane, and another [side] MN is mobile, still always staying perpendicular to AC.
For if the whole AC is divided into two [equal] parts at O, and one foot of the compass is 
placed at O, and another draws a line with [the aperture of] the interval AO, as long
as it moves along the arc CD, let two rules AZ and MN, intersecting at the point D,  be pushed 
until the ruler MN intersects the cursor FY at point E, positioned on the line AC.
By this the same means AE and AD will be obtained, with a much more convenient
method  than the preceding one. This method, if the Figure was much larger, [placed]
somewhere on a huge wall, could be of use to architects.

It is not a place for me to give advise about the compass and its legs moveable is a certain
ruler, or about the sharpness of the other leg, which should delicately push rules
at point D, or  finally about the pressure of rules AZ and MN against the compass. 
This  indeed someone ingenious will attempt to achieve by the weight of rulers,
or by some spring, or by the push of one finger, or will devise something much better than 
that.
\pend
\end{Rightside}

\Columns

\end{pairs}


 \begin{center}
     \includegraphics[scale=0.1]{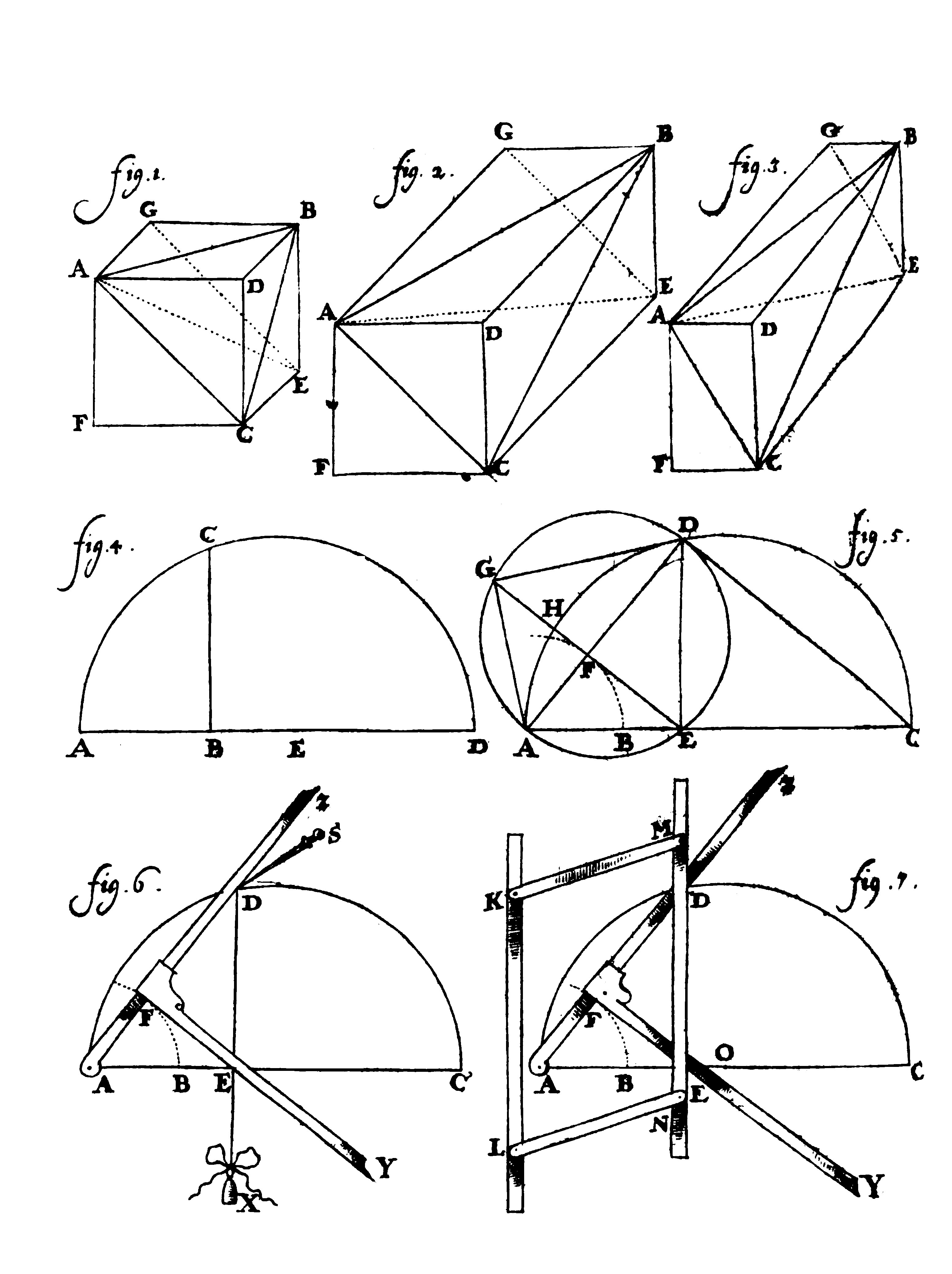}
   \end{center}
 \begin{center}
Figures 1 -- 7.
   \end{center}

\newpage
\section*{Appendix - list of propositions from \emph{Elements} mentioned in the 
text of \emph{Solutio theorematum}}
 
\begin{flushleft}
 \emph{The first number indicates
proposition number, the second one is the book number.
Latin text from \cite{Euclid-latin}, English translation of propositions from \cite{Euclid}.
Pappus generalization of 47.1 and Clavius' scholium for Prop. 31.3 translated by H.F. }
\end{flushleft}

\begin{pairs}
\begin{Leftside}
\selectlanguage{latin}
\pstart
\noindent{\bf Prop. 4.1}\\
Si duo triangula duo latera duobus lateribus aequalia
habeant, alterum alteri; habeant autem et angulum
angulo aequalem, qui aequalibus rectis lineis continentur:
et basim basi aequalem habebunt; et triangulum triangulo
aequale erit; et reliqui anguli reliquis angulis aequales, alter alteri, quibus aequalia latera subtenduntur.

\noindent{\bf Prop. 47.1}\\
In rectangulis triangulis, quod a latere rectum angulum subtendere describitur, quadratum aequale est quadratis quae a lateribus rectum angulum continentibus describuntur.
\pend

\pstart
\noindent{\bf Prop. 12.2}\\
In obtusangulis triangulis quadratum ex latere obtusum angulum subtendere, majus est quam quadrata ex lateribus obtusum angulum continentibus, rectangulo contento bis ab uno laterum quae sunt circa obtusum angulum in quod productum perpendicularis cadit, et recta linea intercepta exterius a perpendiculari ad angulum obtusum.

\pend

\pstart
\noindent{\bf Prop. 13.2}\\
In omni triangulo, quadratum ex latere acutum angulum subtendere, minus est quam quadrata ex lateribus angulum illum continentibus, rectangulo contento bis ab uno laterum quae sunt circa acutum angulum, in quod productum perpendicularis cadit, et recta linea intercepta a perpendiculari ad angulum acutum.
\pend

\pstart
\noindent{\bf Prop. 3.3}\\ 
Si in circulo recta linea per centrum ducta, rectam lineam quandam non ductam per centrum bifariam secet; et ad angulos rectos ipsam secabit quod si ad angulos rectos ipsam secet, et bifariam secabit. 
\pend

\pstart
\noindent{\bf Prop. 8.6}\\
Si in triangulo rectangulo, ab angulo recto ad basim perpendicularis ducatur; quae ad perpendicularem sunt triangula, et toti, et inter se sunt similia.
\pend

\pstart
\noindent{\bf Prop. 31.6}\\
In triangulis rectangulis figura rectilinea quae sit a latere rectum angulum subtendere, aequale est eis quae a lateribus rectum angulum continentibus sunt, similibus et similiter descriptis.
\pend

\pstart
\noindent{\bf Prop. 19.7}\\
Si quatuor numeri proportionales fuerint, qui ex primo, \& quarto fit, numerus,
aequalis erit ei, qui ex secundo, \& tertio fit,
numero. Et si, qui ex primo, \& quarto fit,
numerus, aequalis fuerit ei, qui ex secundo,
\& tertio fit, numero; ipsi quatuor numeri proportionales erunt.

\noindent{\bf Prop. 20.7}\\
Si tres numeri proportionales fuerint; 
qui sub extremis continetur, aequalis est ei,
qui a medio efficitur: Et si, qui sub extremis
 continetur, aequalis fuerit ei, qui a medio describitur; ipsi tres numeri
proportionales erunt.

\noindent{\bf Prop. 4.11}\\
Si recta linea duabus rectis lineis se invicem secantibus in communi sectione ad rectos angulos insistat, etiam ducto per ipsas plano ad rectos angulos erit.
\pend

\pstart
\noindent{\bf Prop. 7.12}\\
Omne prisma triangulem habens basim, dividitur in tres pyramides aequales inter se, quae triangulares bases habent.

\noindent{\bf Pappus generalization of 47.1, as given by Clavius}\\
In omni triangulo, parallelogramma quaecunque super duobus lateribus descripta, aequalia sunt parallelogrammo super reliquo latere constituto, cuius alterum latus aequale sit, \& parallelum  rectae ductae ab angulo, quae duo illa latera comprehendunt, ad punctum, in quo conveniunt latera parallelogrammorum lateribus trianguli opposita, si ad partes anguli dicti producantur. 
\pend

\pstart
\noindent{\bf Scholium Clavii ad 31.3}\\
Manifestum quoque est conversum huius theorematis. Hoc est, segmentum circuli, in quo angulus constitutus est rectus, semicirculus est. Nam si esset maius, angulus in eo foret acutus; si minus, obtusus.
\pend

\pstart
\noindent{\bf Corollarium ad Prop. 8.6}\\
Ex hoc manifestum est, perpendicularem quae in rectangulo triangulo ab angulo recto in basin demittitur, esse mediam proportionalem inter duo basis segmenta.
\pend
\end{Leftside}

\begin{Rightside}
\selectlanguage{english}
\pstart
\noindent{\bf Prop. 4.1}\\
If two triangles have the two sides equal to two sides
respectively, and have the angles contained by the equal straight
lines equal, they will also have the base equal to the base, the
triangle will be equal to the triangle, and the remaining angles
will be equal to the remaining angles respectively, namely those
which the equal sides subtend.

\noindent{\bf Prop. 47.1}\\
In right-angled triangles the square on the side subtending the right angle
is equal to the squares on the sides containing the right angle.
\pend

\pstart
\noindent{\bf Prop. 12.2}\\
In obtuse-angled triangles the square on the side subtending the obtuse
angle is greater than the squares on the sides containing the obtuse angle by twice the
rectangle contained  by one of the sides about the obtuse angle, namely that on which the
perpendicular falls, and the straight line cut off outside by the perpendicular towards
the obtuse angle.
\pend

\pstart
\noindent{\bf Prop. 13.2}\\
In acute-angled triangles the square on the side subtending the acute
angle is less than the squares on the sides containing the acute angle by twice the
rectangle contained  by one of the sides about the acute angle, namely that on which the
perpendicular falls, and the straight line cut off within by the perpendicular towards
the acute angle.
\pend

\pstart
\noindent{\bf Prop. 3.3}\\ 
If in a circle a straight line through the centre bisects a straight line not through 
the centre, it also cuts it at right angles; and if it cut at right angles, it also bisects it.
\pend

\pstart
\noindent{\bf Prop. 8.6}\\
If in a right-angled triangle a perpendicular be drawn from the right angle to the base, the triangles adjoining the perpendicular are similar both to the whole and to one another.
\pend

\pstart
\noindent{\bf Prop. 31.6}\\
In right-angled triangles the figure of the side subtending the right angle is equal to the similar and similarly described
figures on the sides containing the right angle.
\pend

\pstart
\noindent{\bf Prop. 19.7}\\
If four numbers be proportional, the number produced from the first
and fourth will be equal to the number produced from the second and the third;
and, if the number produced from the first and fourth be equal to that produced 
from the second and third, the four numbers will be proportional.

\noindent{\bf Prop. 20.7} (stated in the commentary to prop. 19.7 in \cite{Euclid}) \\
If three numbers be proportional, the product of the extremes
is equal to the square of the mean, and conversely.

\noindent{\bf Prop. 4.11}\\
If a straight line be set up at right angles to two straight
lines which cut one another, at their common point of section,
it will also be at right angles to the plane through them.
\pend

\pstart
\noindent{\bf Prop. 7.12}\\
Any prism which has a triangular base is divided into three
pyramids equal to one another which have triangular bases.

\noindent{\bf Pappus generalization of 47.1, as given by Clavius}\\
In every triangle, any parallelograms built on two sides are equal to the parallelogram built on the remaining side,
whose other side is equal and parallel to the straight line drawn from the angle made by the two sides of the triangle to
the point of intersection of extensions of the sides of parallelograms opposite to the sides of the triangle. 
\pend

\pstart
\noindent{\bf Clavius' scholium for Prop. 31.3}\\
It is also clear that the converse of this theorem is true. That is, a segment of a circle, in which a right angle is constituted,
is a semicircle. For is it was greater, the angle in it would be acute, if lesser, it would be obtuse.
\pend

\pstart
\noindent{\bf Corollary to Prop. 8.6}\\
From this is clear that, if a right angled triangle a perpendicular be drawn from the right angle to the base, the straight line so drawn is a mean proportional between the segments of the base.
\pend
\end{Rightside}

\Columns
\end{pairs}




\begin{thebibliography}{1}

\bibitem{KochanskiSolutio}
Adam~Adamandy Kocha{\'n}ski.
\newblock ``Solutio Theorematum Ab illustri Viro in Actis hujus Anni Mense Januario,
pag. 28. propositorum'',
\newblock {\em Acta Eruditorum}, Lipsiae 1682, pp. 230--236.

\bibitem{Euclid}
Euclid.
\newblock {\em The thirteen books of the Elements}, translated with introduction 
and commentary by Sir Thomas L. Heath,	
\newblock Dover, New York, 1956.

\bibitem{Euclid-latin}
Euclides, Christophorus Clavius.
\newblock {\em Euclidis Elementorum libri XV},
\newblock Romae apud Vincentium Accoltum, 1574.

\end{thebibliography}
\end{document}